\documentclass[11pt]{article}

\usepackage[T1]{fontenc}
\usepackage[latin1]{inputenc}
\usepackage{empheq,amsthm,amssymb}
\usepackage{mathrsfs}
\usepackage{xcolor}
\usepackage{hyperref}
\usepackage[margin=2.55cm]{geometry}
\usepackage{authblk}

\newtheorem{theo}{Theorem}[section]
\newtheorem{prop}[theo]{Proposition}
\newtheorem{cor}[theo]{Corollary}
\newtheorem{lem}[theo]{Lemma}

\numberwithin{equation}{section}

\def\tr{\mbox{\rm Tr}}

\newcommand{\Aa}{ {\cal A }}
\newcommand{\Ba}{ {\cal B }}

\newcommand{\La}{ {\cal L }}

\newcommand{\Ka}{ {\cal K }}
\newcommand{\Ya}{ {\cal Y }}
\newcommand{\Qa}{ {\cal Q }}
\newcommand{\Ua}{ {\cal U }}
\newcommand{\Ga}{ {\cal G }}
\newcommand{\Xa}{ {\cal X }}
\newcommand{\Ha}{ {\cal H }}
\newcommand{\Wa}{ {\cal W }}

\def \PP{\mathbb{P}}
\def \RR{\mathbb{R}}
\def \EE{\mathbb{E}}
\def \QQ{\mathbb{Q}}

\begin{document}

\title{Backward Nonlinear Smoothing Diffusions}
\author[$1$]{Brian D.O. Anderson\thanks{B.D.O. Anderson was supported by the Australian
Research Council (ARC) via grant DP160104500 and grant DP190100887; and by Data61-CSIRO.}}
\author[$2$]{Adrian N. Bishop}
\author[$3$]{Pierre Del Moral\thanks{P. Del Moral was supported in part by the Chair Stress Test, RISK Management and Financial Steering, led by the French Ecole Polytechnique and its Foundation and sponsored by BNP Paribas.}}
\author[$4$]{Camille Palmier}
\affil[$1$]{{\small Research School of Electrical, Energy and
Material Engineering, Australian National University, Canberra, Australia;
and also with Hangzhou Dianzi University, China, and Data61-CSIRO in
Canberra, Australia.}}
\affil[$2$]{{\small CSIRO \& University of Technology Sydney (UTS), Australia}}
\affil[$3$]{{\small INRIA, Bordeaux Research Center, Talence  \& CMAP, Polytechnique Palaiseau, France}}
\affil[$4$]{{\small  Institut de Math\'ematiques de Bordeaux (IMB), Bordeaux University  \& ONERA Palaiseau, France}}
\date{}

\maketitle

\begin{abstract}
We present a backward diffusion flow (i.e. a backward-in-time stochastic differential equation) whose marginal distribution at any (earlier) time is equal to the smoothing distribution when the terminal state (at a latter time) is distributed according to the filtering distribution. This is a novel interpretation of the smoothing solution in terms of a nonlinear diffusion (stochastic) flow. This solution contrasts with, and complements, the (backward) deterministic flow of probability distributions (viz. a type of Kushner smoothing equation) studied in a number of prior works. A number of corollaries of our main result are given including a derivation of the time-reversal of a stochastic differential equation, and an immediate derivation of the classical Rauch-Tung-Striebel smoothing equations in the linear setting.
\end{abstract}

\medskip\noindent\textbf{\small{}Keywords}{\small{}: Nonlinear filtering and smoothing; Kalman-Bucy filter; Rauch-Tung-Striebel smoother; particle filtering and smoothing; diffusion equations; stochastic semigroups; backward stochastic integration; backward It\^o-Ventzell formula; time-reversed stochastic differential equations; Zakai and Kushner-Stratonovich equations.}{\small \par \bigskip}

\noindent\textbf{\small{}Mathematics Subject Classification}{\small{}: 60G35; 62M20; 93E11; 93E14; 60J60.}{\small \par}

\section{Introduction}
Let $(W_t,V_t)\in (\RR^{p}\times\RR^q)$ be a $(p+q)$-dimensional Brownian motion for finite $p,q\geq 1$. Consider a signal-observation model $(X_t,Y_t)\in (\RR^{m}\times\RR^n)$ given by the It\^o stochastic differential equation:
\begin{equation}\label{def-model}
\left\{
\begin{array}{rcl}
dX_t&=&a_t(X_t)\,dt+\sigma_t(X_t)\, dW_t\\
dY_t&=&b_t(X_t)\, dt+\varsigma_t\, dV_t
\end{array}
\right.
\end{equation}
for some measurable functions $\varsigma_t,a_t(x),\sigma_t(x),b_t(x)$ with appropriate dimensions. We set $Y_0=0$ and let $X_0$
be an initial random variable with absolute moments of any order. We let $\alpha_t(x):=\sigma_t(x)\sigma^{\prime}_t(x)$, and $\beta_t:=\varsigma_t\varsigma_t^{\prime}$, where $A^{\prime}$ denotes the transpose of some matrix $A$.

To avoid unnecessary technical details, we assume $\beta_t\geq \epsilon\,I$, for some $\epsilon>0$, where $I$ denotes the identity matrix. We also assume the drift and sensor functions $(a_u(x),b_u(x))$, as well as the diffusion matrix $\sigma_u(x)$, are smooth w.r.t. $(u,x)$ and they have uniformly bounded derivatives w.r.t $x$ of all order on $(u,x)\in [s,t]\times\RR^m$, for any $s\leq t$. 

These technical conditions ensure that the above stochastic differential equation (\ref{def-model}) has a global solution $(X_t,Y_t)$ in the sense of It\^o. In addition, $(X_t,Y_t)$ as well as the sensor function $b_t(X_t)$ have absolute moments of any order.  The stochastic flow associated with the signal is also smooth w.r.t. its initial condition, and its derivatives have absolute moments of any order. 

The filtering problem then consists of computing the conditional distribution $\pi_t$ of the random signal states $X_t$ of the signal 
given the sigma-field $\Ya_t=\sigma(Y_s,~s\leq t)$ generated by the observations. The smoothing problem is to compute the conditional distribution $\pi_{t,s}$ of the random signal states $X_s$ given $\Ya_t$, with $ t\geq s$. With this notation, we have $\pi_{t,t}=\pi_t$. 

The filtering and smoothing problems have been studied extensively, and the literature on this topic is too broad to survey in detail here; and a review of this type is beyond the rather narrow scope of our contribution. We may point to the general texts \cite{Meditch73,bainCrisan} for broad coverage of these problems. 

We do note some rather seminal early literature in the linear setting \cite{bryson,rauch,kailath,rutkowki} and the nonlinear setting \cite{bryson,leondes,anderson-2,andersonRhodes}. The first work on the smoothing topic is the maximum likelihood solution in \cite{bryson} in both the linear and nonlinear setting. The study of \cite{rauch} more formally confirms the linear result in \cite{bryson} and also provides a simpler formulation for the mean and covariance of the smoothing distribution. In the nonlinear setting, the work of \cite{anderson-2,leondes} introduces an analogue of a type of Kushner-Stratonovich equation (see \cite{bainCrisan} for this equation in the filtering context) for the smoothing problem. More specifically, \cite{anderson-2,leondes} propose a deterministic partial differential equation that describes the flow of the smoothing distribution in terms of a backward flow and the standard filtering distribution which acts as the boundary condition (the latter follows from the classical Kushner-Stratonovich equation).

In Section~\ref{main-sec} we state the main contribution of this work. Our main result asserts a backward diffusion flow (i.e. a backward stochastic differential equation) whose marginal distribution at any time $0\leq s\leq t$ is equal to the smoothing distribution $\pi_{t,s}$ when the terminal state is distributed according to the filtering distribution $\pi_t$. 

This is a novel interpretation of the smoothing solution in terms of a nonlinear diffusion (stochastic) flow (in the spirit of McKean-Vlasov-type processes). This solution contrasts with, and complements, say, the (backward) deterministic flow of probability distributions (viz. a type of Kushner smoothing equation) in \cite{anderson-2,leondes}. We also provide a number of corollaries of our main result in Section \ref{illustrations} including an immediate derivation of the Rauch-Tung-Striebel smoothing equations \cite{rauch} in the linear setting.

A number of auxiliary contributions are set forth in order to prove our main contribution to the smoothing problem. As is typical, (e.g. see \cite{bryson,rauch,kailath,rutkowki,leondes,anderson-2,andersonRhodes,pardoux}), our smoothing solution requires the formulation of a related filtering problem. In Section~\ref{nl-filtering} we present a brief review of the Kallianpur-Striebel formula. We then provide a novel and more direct approach to deriving weak-versions of the Zakai and the Kushner-Stratonovich equations in Sections \ref{unnormalized-sg} and \ref{normalized-sg} respectively. We also consider the backward versions of these equations in Section~\ref{backward-stoch-sg}.

Our approach to the filtering equations in this article combines forward and backward It\^o formulas for stochastic transport semigroups with a recent backward version of the It\^o-Ventzell formula presented in~\cite{dp-ss}. This semigroup methodology can be seen as an extension, to the Zakai and Kushner-Stratonovich equations, of the forward-backward stochastic analysis of diffusion flows developed in~\cite{daprato-2,daprato-3,dp-ss,Kunita-82,kunita-b-82}.
 
Our direct semigroup approach to the forward/backward filtering equations in this work contrasts with classical stochastic partial differential methods and functional analysis in Sobolev spaces; see e.g. the seminal works by Pardoux~\cite{pardoux-0,pardoux-2,pardoux}, as well as Krylov and Rozovskii~\cite{krylov-0,krylov-1}. Related reverse time diffusions and filtered and smoothed densities are also developed in~\cite{anderson,andersonRhodes} using discrete time approximation techniques, without a detailed discussion on the existence of these densities. We present a number of auxiliary results in this direction throughout Section~\ref{nl-filtering} which are utilised in the proof of our main smoothing result in Section \ref{backward-sec}.

\subsection{Some preliminary notation}

This subsection presents some notation needed from the onset. 

The signal and the observation defined in (\ref{def-model}) are column vectors. Unless otherwise stated, we use the letters $f$ and $g$ to denote bounded scalar measurable test functions on some measurable space.

We denote by $\nabla f$ the column gradient whenever $f$ is a differentiable function on some Euclidian space, and by $\nabla^2 f$ the Hessian matrix whenever it is twice differentiable.

With $f:\mathbb{R}^m\rightarrow\mathbb{R}$, we let $\mbox{\rm div}_{\alpha_t}(f)$ be the $\alpha_t$-divergence $m$-column vector operator with $j$-th entry given by the formula
\begin{equation*}
\mbox{\rm div}_{\alpha_t}(f)(x)^j:=~\sum_{1\leq i\leq m}~\partial_{x_i}\left(\alpha^{i,j}_t(x)~f(x)\right)
\end{equation*}

The generator $L_t$ of the signal $X_t$ is also given by the second order differential operator
\begin{equation*}
  L_t(f)(x) :=\nabla f(x)^{\prime} b_t(x)+\frac{1}{2}~\tr\left(\nabla^2f(x)\,\alpha_t(x)\right)\quad\mbox{\rm with the trace operator $\tr(\cdot)$.}
\end{equation*}
Here and throughout, and without further mention, we assume that functions $f$ acted on by a second-order differential generator are in addition twice differentiable with bounded derivatives.

For a measure $\mu$ and test function $f$ of compatible dimension we write,
$$
	\mu(f) := \int \mu(dx)\, f(x)
$$ 
An integral operator $\Qa(x,dz)$ acts on the right on scalar test functions $f$; and on the left on measures $\mu$ according to the formulae,
$$
\Qa(f)(x):=\int \Qa(x,dz)~f(z)\quad\mbox{\rm and}\quad (\mu\Qa)(dz):=\int\mu(dx)~\Qa(x,dz)
$$
We extend this operator to an integral operator on matrix functions $h(x)=(h_{i,j}(x))_{i,j}$ 
by setting,
$$
\Qa(h)(x)_{i,j}=\Qa(h_{i,j})(x)
$$

\section{Main Result}\label{main-sec}

In further development of this article we assume for any $t>0$ the conditional distribution $\pi_t$ has a positive density $p_t:=d\pi_t/d\lambda$ w.r.t. the Lebesgue measure $\lambda$ on $\RR^m$. In addition, $p_u(x)$ and its derivative $\nabla p_u(x)$ are uniformly bounded w.r.t. $(u,x)\in [s,t]\times\RR^m$, for any given $s>0$, almost surely w.r.t. the distribution of the observation process. A more detailed discussion on these regularity conditions is provided in Section~\ref{basics}.

The main result of the article takes the following form:

\begin{theo}\label{theo-ref}
For any $ t\geq u\geq s$ we have the transport equation
\begin{equation}\label{theo-eq-pi-K}
\pi_{t,s}(dx)=(\pi_{t,u}\,\Ka_{u,s})(dx):=\int~\pi_{t,u}(dz)~\Ka_{u,s}(z,dx)
\end{equation}
where $\Ka_{u,s}$ denotes the Markov semigroup of the backward diffusion flow,
\begin{equation}\label{theo-eq}
d\Xa_{u,s}(x)=-\left(
\left(p_s(\Xa_{u,s}(x))^{-1}\mbox{\rm div}_{\alpha_s}(p_s)(\Xa_{u,s}(x))-a_s(\Xa_{u,s}(x))\right)\,ds+
\sigma_s(\Xa_{u,s}(x))\,d\Wa_s\right)
\end{equation}
with the boundary condition $\Xa_{u,u}(x)=x$, and where $\Wa_t\in \RR^{p}$ denotes a $p$-dimensional Brownian motion independent of the observations. 
\end{theo}

The proof of the above theorem is provided in Section~\ref{theo-ref-proof}. The backward stochastic differential equation (\ref{theo-eq}) should be read as shorthand for the backward It\^o integration formula,
\begin{equation}\label{theo-eq-int}
\Xa_{t,s}(x)=x+\int_s^t\left(p_u(\Xa_{t,u}(x))^{-1}\mbox{\rm div}_{\alpha_u}(p_u)(\Xa_{t,u}(x)) -a_u(\Xa_{t,u}(x))
\right)\,du+\int_s^t
\sigma_u(\Xa_{t,u}(x))\,d\Wa_u
\end{equation}
with terminal condition $\Xa_{t,t}(x)=x$. The right-most term in the above formula is an It\^o backward stochastic integral such that for any terminal time $t$ this process is a square integrable backward martingale w.r.t. the variable $s\in [0,t]$. 

Formally, we may slice the time interval $[s,t]_h:=\{u_0,\ldots,u_{n-1}\}$ via some time mesh $u_{i+1}=u_i+h$ from $u_0=s$ to $u_{n}=t$ and with time step $h>0$. In this notation, according to the backward equation (\ref{theo-eq}), or (\ref{theo-eq-int}), we compute $\Xa_{t,u-h}(x)$ from $\Xa_{t,u}(x)$ using the formula
\begin{equation}\label{fb-eq-discrete}
\begin{array}{l}
\Xa_{t,u-h}-\Xa_{t,u}~\simeq~
\left(p_u(\Xa_{t,u})^{-1}\mbox{\rm div}_{\alpha_u}(p_u)(\Xa_{t,u})
-a_u(\Xa_{t,u})\right)h+\sigma_u(\Xa_{t,u})(\Wa_u-\Wa_{u-h})
\end{array}
\end{equation}

We may provide some comments on the above theorem. By construction, given the observations and for any given $x\in\RR^m$ and $t\geq s$, the probability $\Ka_{t,s}(x,dz)$ introduced in (\ref{theo-eq-pi-K}) coincides with the distribution of the random state $\Xa_{t,s}(x)$.
In addition, for any $t\geq u\geq s$ we have the integral and stochastic semigroup properties,
\begin{equation}
\Ka_{t,s}(x_2,dx_0):=\int~\Ka_{t,u}(x_2,dx_1)~\Ka_{u,s}(x_1,dx_0)
\end{equation}
and
\begin{equation}
\quad\Xa_{t,s}=\Xa_{u,s}\circ\Xa_{t,u}
\end{equation}
where $\Xa_{u,s}\circ\Xa_{t,u}$ denotes the composition of the mappings $\Xa_{u,s}$ and $\Xa_{t,u}$.

If we let $\Xa_t$ be a random variable with distribution $\pi_t$, for some $t\geq 0$. According to (\ref{theo-eq-pi-K}) the  random state  $\Xa_{t,s}(\Xa_t)$ of the process (\ref{theo-eq}) at any given $s\in [0,t]$, is distributed according to $\pi_{t,s}=\pi_{t}\,\Ka_{t,s}$. In words, the backward process $\Xa_{t,s}(\Xa_t)$ is distributed according to the smoothing distribution $\pi_{t,s}$ for any $s\leq t$ whenever the terminal condition $\Xa_{t,t}(\Xa_t)=\Xa_t$ is distributed according to the filtering distribution $\pi_t$. In this sense, (\ref{theo-eq}) constitutes a backward nonlinear smoothing diffusion. A forward diffusion flow that has a marginal distribution at any time equal to the filtering distribution is considered in \cite{yang2013TAC,Yang2016}.

More generally, we have the backward It\^o formula
\begin{equation}\label{ref-proof}
df(\Xa_{t,s}(x))=-\La_{s,\pi_s}(f)(\Xa_{t,s}(x))~ds-\nabla f(\Xa_{t,s}(x))^{\prime}~\sigma_s(\Xa_{t,s}(x))~d\Wa_s
\end{equation}
with the second order differential operator
\begin{equation}
 \La_{s,\pi_s}(f)=\sum_{1\leq j\leq m}\left(-\,a^j_s+\frac{1}{p_s}~\mbox{\rm div}_{\alpha_s}(p_s)^j\right)~
 \partial_{x_j}f+\frac{1}{2}\sum_{1\leq i,j\leq m}\alpha^{i,j}_s~ \partial_{x_ix_j}f
\end{equation}
Equivalently, we have the backward martingale decomposition
\begin{equation}
f(\Xa_{t,s}(x))-f(x)-\int_s^t\La_{u,\pi_u}(f)(\Xa_{t,u}(x))~du=\int_s^t\nabla f(\Xa_{t,u}(x))^{\prime}~\sigma_u(\Xa_{t,u}(x))~d\Wa_u
\end{equation}
This yields the backward evolution equations,
\begin{equation}
\partial_s\Ka_{t,s}(f)(x)=
-\Ka_{t,s}(\La_{s,\pi_s}(f))(x)
\end{equation}
and
\begin{equation}\label{fb-eq-leondes}
\partial_{s}\pi_{t,s}(f)=-\pi_{t,s}(\La_{s,\pi_s}(f))
\end{equation}
with terminal conditions $\Ka_{t,t}(f)=f$ and $\pi_{t,t}=\pi_t$. The formula (\ref{fb-eq-leondes}) coincides with the conditional Fokker-Planck equation in~\cite{leondes}, and further developed in~\cite{anderson-2}.

For further discussion on general backward integration of stochastic flows see~\cite{daprato-2}; see also the appendix of~\cite{mp-var-18} in the context of nonlinear diffusions, and~\cite{pardoux} in the context of nonlinear filtering, and~\cite{dp-ss} on forward-backward perturbation analysis of stochastic flows. Note there is no issue with adaption of the backward process in the sense studied in \cite{pardoux-peng} since we rely only on the independent backward Brownian motion in (\ref{theo-eq}). The ``backward diffusion'' in (\ref{theo-eq}) is backward in the sense of a time reversed stochastic differential equation as in \cite{anderson,haussmann,millet}.

\subsection{Some corollaries}\label{illustrations}

We end this introduction with some direct consequences of the above theorem.

Note when $b_t=0$ the measure $\pi_t$ coincides with the distribution of the random state $X_t$ of the signal. In this context, $\Xa_{t,s}(X_t)$ reduces to the time reversal of the signal starting at $\Xa_{t,t}(X_t)=X_t$ at the terminal time $t$. Using Theorem~\ref{theo-ref} we recover the fact that the time reversal process of the signal is itself a Markov diffusion~\cite{anderson,haussmann,millet}. More precisely, we have the corollary:

\begin{cor}[Anderson~\cite{anderson}]
Assume that $b_t=0$. For any time horizon $t\geq 0$, the process
$\mathfrak{X}^t_s:=X_{t-s}$ with $s\in [0,t]$ is a Markov process with generator
\begin{equation}
\mathfrak{L}_{s}^t(f)=
\sum_{1\leq j\leq m}\left(\frac{1}{p_{t-s}}~\mbox{\rm div}_{\alpha_{t-s}}(p_{t-s})^j-\,a^j_{t-s}\right)~
 \partial_{x_j}f+\frac{1}{2}\sum_{1\leq i,j\leq m}\alpha^{i,j}_{t-s}~ \partial_{x_ix_j}f
\end{equation}
\end{cor}

We consider now linear-Gaussian filtering/smoothing models with,
\begin{equation}\label{lg-model}
a_t(x)=A_t\,x, \qquad b_t(x)=B_t\,x \quad\mbox{\rm and homogeneous diffusion matrix}\quad\sigma_t(x)=\Sigma_t
\end{equation}
for some matrices $(A_t,B_t,\Sigma_t)$ with appropriate dimensions. Whenever $X_0$ is a Gaussian random variable with mean $\widehat{X}_0$ and covariance matrix $R_0$, the optimal filter $\pi_t$ is a Gaussian distribution with mean $\widehat{X}_t$ and covariance matrix $R_t$ satisfying the Kalman-Bucy and the Riccati equations 
\begin{equation}\label{ref-kalman-riccati}
\left\{
\begin{array}{rcl}
d\widehat{X}_t&=&A_t\widehat{X}_t~dt+R_tB_t^{\prime}\beta_t^{-1}\left(dY_t-B_t \widehat{X}_t~dt\right)\\
\partial_tR_t&=&A_tR_t+R_tA^{\prime}_t+\alpha_t-R_tB_t^{\prime}\beta_t^{-1}B_t R_t
\end{array}
\right.
\end{equation}
In this context, we also have that
\begin{equation}\label{ref-growth}
-p_s(x)^{-1}\mbox{\rm div}_{\alpha_s}(p_s)(x)=\alpha_sR^{-1}_s(x-\widehat{X}_s)
\end{equation}
This yields the following corollary:
\begin{cor}
For linear Gaussian filtering models (\ref{lg-model}), the diffusion flow $\Xa_{t,s}(x)$ satisfies the backward evolution equation
\begin{equation}\label{ref-GL}
d\Xa_{t,s}(x)=-\left(
\left(-A_s~\Xa_{t,s}(x)-\alpha_sR^{-1}_s(\Xa_{t,s}(x)-\widehat{X}_s)\right)\,ds+
\Sigma_s~d\Wa_s\right)
\end{equation}
with the boundary condition $\Xa_{t,t}(x)=x$. 
\end{cor}

Replacing $x$ in (\ref{ref-GL}) by a random variable $\Xa_{t}$ with distribution $\pi_t$ for any $t\geq s$ we have that $\Xa_{t,s}(\Xa_{t})$ has distribution $\pi_{t,s}$. In addition, since the process is linear the distribution $\pi_{t,s}$ is Gaussian with mean $\widehat{X}_{t,s}$ and covariance matrix $R_{t,s}$. The discrete time version of (\ref{ref-GL}) can be found in Section 9.9.6 in~\cite{penev}.

 Now taking expectations we readily deduce the rather well-known Rauch-Tung-Striebel smoothing equations~\cite{rauch}, simplifying the innovation techniques and the sophisticated approximation theory developed in~\cite{kailath,leondes,rutkowki}, or the formal variational approaches and maximum likelihood techniques presented in the pioneering articles~\cite{bryson,rauch}. 

\begin{cor}[Rauch-Tung-Striebel~\cite{rauch}]
For any $t\geq s$, the parameters  $(\widehat{X}_{t,s},R_{t,s})$ satisfy
the backward evolution equations
\begin{equation}
\left\{
\begin{array}{rcl}
\partial_s\widehat{X}_{t,s}&=&A_s\widehat{X}_{t,s}+\alpha_sP^{-1}_s(\widehat{X}_{t,s}-\widehat{X}_s)\\
\partial_sR_{t,s}&=&(A_s+\alpha_sR^{-1}_s)R_{t,s}+R_{t,s}(A_s+\alpha_sR^{-1}_s)^{\prime}-\alpha_s
\end{array}
\right.
\end{equation}
with terminal conditions $(\widehat{X}_{t,t},R_{t,t})=(\widehat{X}_{t},R_{t})$.
\end{cor}

\subsection{Comments on our regularity conditions}\label{basics}

We end this section with some comments on the regularity conditions discussed at the beginning of Section~\ref{main-sec}. These conditions are clearly met for linear Gaussian filtering models (see e.g. (\ref{ref-kalman-riccati}) and (\ref{ref-growth})). They are also met for nonlinear models as soon as the signal satisfies a classical controllability-type condition. 
 
Note firstly, whenever the signal is uniformly elliptic, in the sense that $\alpha_t(x)=\sigma_t(x)\sigma^{\prime}_t(x)\geq \delta\,I$ for some $\delta>0$, then it is well known that $X_t$ has a smooth positive density w.r.t. the Lebesgue measure on $\RR^m$. Nevertheless in many important applications this ellipticity condition is not satisfied. The  parabolic H\"ormander condition for time varying models~\cite{cattiaux,hormander} is a weaker condition. For linear-Gaussian filtering problems,  this condition reduces to the usual controllability condition. Indeed, if we replace the Brownian motions
$W_t$ by some arbitrary smooth control functions, all states are accessible from one to another, as soon as the Lie algebra generated by the 
controlled vector fields is of full rank. This result is also called the Chow-Rashevskii theorem~\cite{chow,rashevskii}.  Under this H\"ormander condition, the H\"ormander theorem~\cite{hormander}
ensures that the signal states have a smooth density w.r.t. the Lebesgue measure on $\RR^m$. In addition, for any $s<t$ the Markov transition semigroup 
$P_{s,t}$
of the signal has a smooth positive density $(x,z)\mapsto p_{s,t}(x,z)$ w.r.t. the Lebesgue measure $\lambda$ on $\RR^m$. In addition, the integral operator $P_{s,t}$ with $s<t$ maps test functions  $f$ into bounded smooth functions  $P_{s,t}(f)$ given by
$$
P_{s,t}(f)(x)=\int~P_{s,t}(x,dz)~f(z)=\int~f(z)~p_{s,t}(x,z)~dz
$$

A natural way to transfer the smoothing properties of $P_{s,t}$ to the optimal filter is to use the following equation
\begin{equation}\label{ref-pi-P}
\pi_t(f)=\pi_0(P_{0,t}(f))+\int_0^t\pi_s(P_{s,t}(f)~(b_s-\pi_s(b_s)))^{\prime}~\beta_s^{-1}~(dY_s-\pi_s(b_s)~ds)
\end{equation} 
given in Theorem 1.1 in~\cite{kunita-71}. Using this formula we readily check that for any $t>0$ the conditional distribution $\pi_t$ has a positive density $p_t$ on $\RR^m$. 
Whenever $\sigma_t(x)$ and $b_t(x)$ are also bounded, Theorem  3.6 in~\cite{michel} (see also Theorem 6.3 in~\cite{Kunita-82}) also ensures that $p_u$ is smooth, and for any $k\geq 1$, any parameters $h>0$ and any time horizon $t>0$ we have
\begin{equation}\label{ref-hyp}
\sup_{h\leq s\leq t}{\,\sup_{x\in \RR^m}{\left(\vert p_s(x)\vert+\Vert\nabla^k p_s(x)\Vert\right)}}<\infty
\end{equation} 
where $\Vert\cdot\Vert$ stands for any (equivalent) norm on $\RR^m$.  
 
 The above estimates are met for linear Gaussian filtering models. Nevertheless,
some caution must be used when considering estimates of the form (\ref{ref-hyp}). Indeed, most of the literature on  stochastic partial differential equations arising in nonlinear filtering, such as the strong formulation of the Zakai and Kushner-Stratonovich equations, assume that the sensor function is uniformly bounded, see e.g.~\cite{Kunita-82,michel,pardoux} and~\cite{veretennikov,zakai}. To the best of our knowledge the extension of the estimate (\ref{ref-hyp}) to more general unbounded sensor functions is still an open important question.

We also note here that the Kallianpur-Striebel formula \cite{kallianpur-1,kallianpur-2} is valid as soon as $\beta_u\geq \epsilon\,I$, for some $\epsilon>0$ and the functions $(a_u(x),b_u(x),\sigma_u(x))$ are smooth with uniformly bounded derivatives w.r.t $x$ of all order on $(u,x)\in [s,t]\times\RR^m$, for any $s\leq t$. Weaker conditions can also be found in the book~\cite{bainCrisan} and the recent article~\cite{cass}.

	Since $X_t$ has continuous paths, for any continuous function $f$ and any $s\leq t$ the random mapping $u\in [s,t]\mapsto f(X_u)$ is almost surely a uniformly bounded function. In addition, $f(X_t)$ is integrable as soon as $f$ has polynomial growth. Up to some classical localization procedure (see e.g. Chapter 7 in~\cite{steele}), these rather weak regularity properties also ensure that the integral semigroups that transport (in time) the filtering measures discussed in Section~\ref{nl-filtering}, as well as their  stochastic partial differential evolution equations, are well defined on any test function with polynomial growth.

\section{Nonlinear filtering equations}\label{nl-filtering}

As is  well known (e.g. see \cite{bryson,rauch,kailath,rutkowki,leondes,anderson-2,andersonRhodes,pardoux}), a solution to the smoothing problem will typically make use of the solution of a related filtering in some way. Consequently, we need to present and develop some related filtering results for proving our main result, Theorem \ref{theo-ref}. This section is largely self-contained but it is vital in the proof, in Section \ref{backward-sec}, of our main result.

The first part of this section presents the classical Kallianpur-Striebel formula which acts as a continuous-time version of Bayes law. In Sections \ref{unnormalized-sg} and \ref{normalized-sg} respectively we present the Zakai, and Kushner-Stratonovich equations for the flow of the conditional filtering distributions (both unnormalised and normalised). These results are rather well known. For further background on these classical ideas, we refer to the pioneering articles by Kallianpur and Striebel~\cite{kallianpur-1,kallianpur-2}, and by  Kushner~\cite{kushner} and Zakai~\cite{zakai}. For more recent discussion on these probabilistic models, we refer to~\cite{cass}, and~\cite{bainCrisan,penev}, and the references therein. In this article, we present a novel and self contained derivation based on stochastic transport semigroups and their forward evolution equations. 

The solution of the Zakai equation is sometimes termed the unnormalized filter.  The 
 semigroup that transports these filtering measures (in time) is discussed in Section~\ref{unnormalized-sg}; and its normalized version in Section~\ref{normalized-sg}. Section \ref{backward-stoch-sg} presents a novel direct approach for deriving the backward evolution of these transport semigroups. Our approach in Section \ref{backward-stoch-sg} combines the backward It\^o formula for stochastic flows with the backward It\^o-Ventzell formula presented in~\cite{dp-ss}.

Now, we introduce some notation/terminology and briefly present the Kallianpur-Striebel formula and the linear semigroup property of unnormalized measures. Let $X_{s,t}(x)$ be the stochastic flow of the signal on the time interval $[s,t]$ and starting at $x$ at time $s$. Let $Z_{s,t}(x)$ be the multiplicative functional
\begin{equation}\label{Z}
 \log{Z_{s,t}(x)}:=\int_s^t b_u(X_{s,u}(x))^{\prime}\beta^{-1}_u~dY_u-\frac{1}{2}\int_s^t\,b_u(X_{s,u}(x))^{\prime}\beta^{-1}_u b_u(X_{s,u}(x))~du
\end{equation}
When $x$ is replaced by $X_s$ we may write $Z_{s,t}$ instead of $Z_{s,t}(X_s)$, and when $s=0$, we may also write $Z_t$ instead of $Z_{0,t}$. With this notation, we have the classical Kallianpur-Striebel formula, 
\begin{equation*}
\pi_t(f)={\gamma_t(f)}/{\gamma_t(1)}\quad\mbox{\rm with}\quad \gamma_t(f):=\EE_0\left(f(X_t)~Z_t\right)
\end{equation*}
Here, $\EE_0(\cdot)$ denotes the expectation operator w.r.t. the signal with a fixed observation process. 
  
The transport semigroup of the unnormalized measures $\gamma_t$ is given for any $s\leq t$ by the formula
\begin{equation}\label{zakai-transport}
\gamma_t=\gamma_sQ_{s,t}
\quad\mbox{\rm
with}\quad
Q_{s,t}(f)(x):=\EE_0\left(f(X_{s,t}(x))~Z_{s,t}(x)\right)
\end{equation}
To check this claim observe that,
$$
Z_t=Z_s\,Z_{s,t}\quad\Longrightarrow \quad \EE_0\left(f(X_t)~Z_t\right)=\EE_0\left(Z_s~ \EE_0(f(X_t)~Z_{s,t}~|~X_s)\right)=\EE_0\left(Z_s~Q_{s,t}(f)(X_s)\right)
$$
Now for any $s\leq u\leq t$ we have
\begin{eqnarray*}
Q_{s,t}(f)(X_s)~=~\EE_0(f(X_t)~Z_{s,t}~|~X_s)&=&\EE_0(Z_{s,u}~\EE\left(f(X_t)~~Z_{u,t}~|~X_u)~|~X_s\right)\\
&=&\EE_0(Z_{s,u}~Q_{u,t}(f)(X_u)~|~X_s)~=~Q_{s,u}(Q_{u,t}(f))(X_u)
\end{eqnarray*}
This yields the integral semigroup formula
$$
Q_{s,t}(x_0,dx_2)=(Q_{s,u}Q_{u,t})(x_0,dx_2):=\int~Q_{s,u}(x_0,dx_1)~Q_{s,u}(x_1,dx_2)
$$
In a more compact form, the semigroup property takes the form
$$
Q_{s,t}=Q_{s,u}Q_{u,t}\quad \mbox{\rm with}\quad Q_{t,t}=I\quad\mbox{\rm where $I$ denotes the identity operator. }
$$
\subsection{Unnormalized stochastic semigroups}\label{unnormalized-sg}
Consider the stochastic transport semigroups $\PP_{s,t}$ and $\QQ_{s,t}$ defined  by
the composition of functions
$$
\PP_{s,t}(f)(x):=(f\circ X_{s,t})(x)~\quad\mbox{\rm and}\quad
\QQ_{s,t}(f)(x):=\PP_{s,t}(f)(x)~Z_{s,t}(x)
$$
Using the semigroup properties of the stochastic flow $X_{s,t}(x)$ for any $s\leq u\leq t$ we check that
$$
\PP_{s,t}(f)(x)=(f\circ X_{s,t})(x)=(f\circ X_{u,t})(X_{s,u}(x))=\PP_{s,u}(\PP_{u,t}(f))(x)
$$
Similarly, we have
$$
\QQ_{s,t}(f)(x)=Z_{s,u}(x)~\left(Z_{u,t}(X_{s,u}(x))~(f\circ X_{s,t})(X_{s,u}(x))\right)=\QQ_{s,u}(\QQ_{u,t}(f))(x)
$$
In a more compact form we have the semigroup properties
$$
\PP_{s,t}=\PP_{s,u}\circ \PP_{u,t}\quad\mbox{\rm and}\quad \QQ_{s,t}=\QQ_{s,u}\circ \QQ_{u,t}\quad \mbox{\rm with}\quad
\PP_{t,t}=I=\QQ_{t,t}
$$
Also observe that
$$
P_{s,t}(f)(x):=\EE_0\left(\PP_{s,t}(f)(x)\right)
\quad\mbox{\rm and}\quad
Q_{s,t}(f)(x):=\EE_0\left(\QQ_{s,t}(f)(x)\right)
 $$
The forward evolution equations of the above semigroups  are described in the next proposition.
\begin{prop}\label{prop-zakai}
For any
 $t\geq s$ we have the forward stochastic evolution equation
 \begin{equation}\label{zakai-ps-ref}
 d\QQ_{s,t}(f)=\QQ_{s,t}(L_t(f))~dt+\QQ_{s,t}(f\,b_t^{\prime})~\beta^{-1}_t~dY_t+\QQ_{s,t}((\nabla f)^{\prime}~\sigma_t)~dW_t
 \end{equation}
  with initial condition $\QQ_{s,s}(f)=f$, when $t=s$. In particular, we have the forward equation
\begin{equation}\label{zakai}
 dQ_{s,t}(f)=Q_{s,t}(L_t(f))\,dt+Q_{s,t}(f\,b_t^{\prime})~\beta_t^{-1}\,dY_t
\end{equation}
with the initial condition $Q_{s,s}(f)=f$, when $t=s$.
\end{prop}
\proof
Assume that the sensor function $b_u(x)$ is uniformly bounded on $[s,t]\times\RR^m$, for any $s\leq t$. In this situation,  the random process $(X_{s,u}(x),Z_{s,u}(x))$ also has uniformly bounded absolute moments of any order on any compact interval $[s,t]$, for any time parameters $s\leq t$.
In this context, we use It\^o formula to check that
$$
dZ_{s,t}(x)=Z_{s,t}(x)~b_t(X_{s,t}(x))^{\prime}\beta^{-1}_t~dY_t$$
as well as
$$
d\PP_{s,t}(f)(x)=\PP_{s,t}(L_t(f))(x)~dt+\PP_{s,t}\left(\nabla f^{\,\prime}~\sigma_t\right)(x)~dW_t
$$
An integration by parts yields
\begin{eqnarray*}
d\QQ_{s,t}(f)(x)&=&Z_{s,t}(x)~d\PP_{s,t}(f)(x)+\PP_{s,t}(f)(x)~dZ_{s,t}(x)\\
&=&\left(L_t(f)(X_{s,t}(x))~Z_{s,t}(x)~dt+Z_{s,t}(x)~f(X_{s,t}(x))~b_t(X_{s,t}(x))^{\prime}\beta^{-1}_t~dY_t\right)\\
&&\hskip3cm+Z_{s,t}(x)~\nabla f(X_{s,t}(x))^{\prime}~\sigma_t(X_{s,t}(x))~dW_t
\end{eqnarray*}
By classical localization principles  of It\^o integrals (see for instance Chapter 7 in~\cite{steele}), the above result is also true for unbounded sensor functions. 
This ends the proof of (\ref{zakai-ps-ref}).
Taking the expectations, we conclude that
\begin{eqnarray*}
d\EE_0(\QQ_{s,t}(f)(x))&=&\EE_0\left(\QQ_{s,t}(L_t(f))(x)\right)~dt+\EE_0\left(\QQ_{s,t}(f~b_t^{\prime})(x)\right)\beta^{-1}_t~dY_t
\end{eqnarray*}
This ends the proof of (\ref{zakai}). The proof of the proposition is completed.\qed

Combining (\ref{zakai-transport}) with Fubini's theorem,  we readily check  the weak form of Zakai equation given
by the formula
\begin{equation}\label{zakai-r}
 d\gamma_t(f)=\gamma_t(L_t(f))\,dt+\gamma_t(f~b_t^{\prime})~\beta_t^{-1}\,dY_t
\end{equation}

 Arguing as in (\ref{ref-pi-P}), we transfer the smoothing properties of $P_{s,t}$  to $Q_{s,t}$ using the perturbation formulae
given for any $s<t$  by
$$
Q_{s,t}(f)=P_{s,t}(f)+\int_s^t~Q_{s,u}\left(P_{u,t}(f)~b_u^{\prime}\right)~\beta_u^{-1}\,dY_u
$$
Arguing as in~\cite{zakai}, the above formula shows that for any $s<t$ the integral operator $Q_{s,t}(x_0,dx_1)$ has a density $x_1\mapsto q_{s,t}(x_0,x_1)$ w.r.t. the Lebesgue measure  on $\RR^m$ given by the integral equation
\begin{equation}\label{perturbation-q-p}
q_{s,t}(x_0,x_1)=p_{s,t}(x_0,x_1)+\int_s^t~\left[\int~q_{s,u}(x_0,z)~p_{u,t}(z,x_1)~b_u^{\prime}(z)~dz\right]~\beta_u^{-1}\,dY_u
\end{equation}

\subsection{Normalized stochastic semigroups}\label{normalized-sg}

 Let $\overline{Z}_{s,t}(x)$ be  the multiplicative functional  defined as $Z_{s,t}(x)$ by replacing in (\ref{Z}) the function 
 $b_u$ and the observation increment $dY_u$ by the centered function $ \overline{b}_u$ and the innovation 
  increment $d\overline{Y}_u$ defined by the formulae
 \begin{equation*}
 \overline{b}_u:=b_u-\pi_u(b_u)\quad \mbox{\rm and }\quad
d \overline{Y}_u:=dY_u-\pi_u(b_u)~du
 \end{equation*}
 Under our assumptions, the random process $\pi_t(b_t)$ is almost surely square integrable on any compact time interval so that the innovation process is well defined.
  Choosing $f=1$ in (\ref{zakai-r}) we check that 
 \begin{equation*}
\log{ \gamma_t(1)}=\int_0^t \pi_u(b_u)^{\prime}\beta^{-1}_u~dY_u-\frac{1}{2}\int_0^t\,\pi_u(b_u)^{\prime}\beta^{-1}_u \pi_u(b_u)~du
 \end{equation*}
 Observe that
 $$
 \pi_sQ_{s,t}(1)=\gamma_t(1)/\gamma_s(1)=\exp{\left(\int_s^t \pi_u(b_u)^{\prime}\beta^{-1}_u~dY_u-\frac{1}{2}\int_s^t\,\pi_u(b_u)^{\prime}\beta^{-1}_u \pi_u(b_u)~du\right)}
 $$
  We also consider the normalized stochastic semigroup
$$
\overline{\QQ}_{s,t}(f)(x):=(f\circ X_{s,t})(x)~\overline{Z}_{s,t}(x)=\PP_{s,t}(f)(x)~\overline{Z}_{s,t}(x)
$$
Arguing as above,  for any $s\leq u\leq t$ we check that
$$
\overline{\QQ}_{s,t}=\overline{\QQ}_{s,u}\circ \overline{\QQ}_{u,t}\quad \mbox{\rm and}\quad\overline{Z}_{s,t}(x)={Z_{s,t}(x)}/{\pi_sQ_{s,t}(1)}
$$

Consider the semigroup
 \begin{equation*}
 \overline{Q}_{s,t}(f)(x):=\EE_0\left(\overline{\QQ}_{s,t}(f)(x)\right)=\EE_0\left(f(X_{s,t}(x))~\overline{Z}_{s,t}(x)\right)={Q_{s,t}(f)(x)}/{\pi_sQ_{s,t}(1)}
 \end{equation*}
 
In this notation, using the same arguments as in the proof of Proposition~\ref{prop-zakai} we have the following forward evolution equations.
 \begin{prop}
 For any given time horizon $s$ and for any
 $t\geq s$ we have the forward stochastic evolution equation
 \begin{equation*}
 d \overline{\QQ}_{s,t}(f)= \overline{\QQ}_{s,t}(L_t(f))~dt+ \overline{\QQ}_{s,t}(f\,\overline{b}_t^{\prime})~\beta^{-1}_t~d\overline{Y}_t+ \overline{\QQ}_{s,t}((\nabla f)^{\prime}~\sigma_t)~dW_t
 \end{equation*}
  with initial condition $ \overline{Q}_{s,s}(f)=f$, when $t=s$. In particular, we have the forward equation
$$
 d \overline{Q}_{s,t}(f)= \overline{Q}_{s,t}(L_t(f))\,dt+ \overline{Q}_{s,t}(f~ \overline{b}_t^{\prime})~\beta_t^{-1}\,d\overline{Y}_t
$$
with the initial condition $\overline{Q}_{s,s}(f)=f$, when $t=s$.
\end{prop}
The above proposition yields the weak form of the Kushner-Stratonovich equation defined  by 
\begin{equation}\label{KS}
d\pi_t(f)=\pi_t(L_t(f))\,dt+ \pi_{t}(f~ \overline{b}_t)^{\prime}~\beta_t^{-1}\,d\overline{Y}_t
\end{equation}
Formally, using the same notation as in (\ref{ref-time-mesh}) we have the forward approximation equation
\begin{equation}\label{ref-time-over-pi}
\pi_{u+h}(f)\simeq \pi_{u}(f)+\pi_{u}(L_u(f))\,h+ \pi_{u}(f~ \overline{b}_u)^{\prime}~\beta_u^{-1}\,\left(\overline{Y}_{u+h}-\overline{Y}_u\right)
\end{equation}

\subsection{Backward evolution equations}\label{backward-stoch-sg}
This section is concerned with the backward evolution equation associated with the unnormalized semigroup $\QQ_{s,t}$ and its normalized version. The main result of this section is the following theorem:

\begin{theo}\label{zakai-Q}
For any twice differentiable function $f$ with bounded derivatives and for any $s\leq t$ we have the backward evolution equation
 \begin{eqnarray}
  d\QQ_{s,t}(f)(x) &=&-\left(\nabla\,\QQ_{s,t}(f)(x)^{\prime}~ a_s(x)+\frac{1}{2}~\tr\left(\nabla^2\,\QQ_{s,t}(f)(x)~ \alpha_s(x)\right)\right)~ds\nonumber\\
&&\hskip1.7cm -\,\QQ_{s,t}(f)(x)~b_s(x)^{\prime}\beta_s^{-1}~dY_s-\nabla\, \QQ_{s,t}(f)(x)^{\prime}~\sigma_s(x)~dW_s  \label{backward-zakai-ps}
 \end{eqnarray}
 with terminal condition $\QQ_{t,t}(f)=f$, when $s=t$. In particular, we have the backward equation
\begin{equation}\label{backward-zakai}
 d Q_{s,t}(f)=-\left(L_s(Q_{s,t}(f))\,ds+Q_{s,t}(f)~b_s^{\prime}~\beta_s^{-1}\,dY_s\right)
\end{equation}
with terminal condition $Q_{t,t}(f)=f$, when $s=t$.
\end{theo}
\proof

We use a direct approach combining
 the backward filtering calculus developed in~\cite{Kunita-82,veretennikov} based on the backward  It\^o calculus developed in~\cite{daprato-2,daprato-3,krylov,kunita-84}, see also the more recent article~\cite{dp-ss} and references therein. 
  
 Consider the stochastic flow $\chi_{s,t}\left(\overline{x}\right)$ starting at $$\chi_{s,s}\left(
\overline{x}\right)=\overline{x}:=\left(
 \begin{array}{c}
 x\\
z
 \end{array}\right)\in (\RR^{m}\times\RR)$$ on the time interval $[s,\infty[$ and given for any $t\geq s$ by
 $$
 \chi_{s,t}\left(
\overline{x}\right):=\left(\begin{array}{c}
 X_{s,t}(x)\\
  Z_{s,t}(x)~z
 \end{array}\right)\in (\RR^{m}\times\RR)
 $$
We
 set
 $$
 \begin{array}{rclcrcl}
 \Ba_t\left(
\overline{x}\right)&:=&\left(
 \begin{array}{c}
 a_t(x)\\
0
 \end{array}\right)&&\Ua_t&:=&\left(
 \begin{array}{c}
W_t\\
Y_t
 \end{array}\right)
\\
 \\
\Lambda_t\left(
\overline{x}\right)&:=&\left(
 \begin{array}{cc}
 \sigma_t(x)&0\\
0&z~b_t(x)^{\prime}\beta_t^{-1}
 \end{array}\right)& \mbox{\rm and }&\Aa_t
\left(
\overline{x}\right)&:=&\Lambda_t\left(
\overline{x}\right)\Lambda_t\left(
\overline{x}\right)^{\prime}
 \end{array} 
 $$
Assume that the sensor function $b_u(x)$ is uniformly bounded on $[s,t]\times\RR^m$, for any $s\leq t$. Then, the process $\left(Z_{s,u}(x), \chi_{s,u}\left(
\overline{x}\right)\right)$ has continuous partial derivatives and also  has uniformly bounded absolute moments of any order on $\left([s,t]\times\RR^m\right)$, for any $s\leq t$. In this situation, we have the forward stochastic evolution equation
 $$
 d \chi_{s,t}\left(
\overline{x}\right)=\Ba_{t}\left(\chi_{s,t}\left(
\overline{x}\right)\right)~dt+\Lambda_t\left(\chi_{s,t}\left(
\overline{x}\right)\right)~d\Ua_t $$
For any twice differentiable function $F$ on $(\RR^{m}\times\RR)$ with bounded derivatives we also have the backward equation
 $$
 \begin{array}{l}
 \displaystyle d(F\circ\chi_{s,t})\left(
\overline{x}\right) \\
\\
 \displaystyle=-\left(\nabla (F\circ\chi_{s,t})\left(
\overline{x}\right)^{\prime}\, \Ba_s\left(
\overline{x}\right)+\frac{1}{2}\,\tr\left(\nabla^2 (F\circ\chi_{s,t})\left(
\overline{x}\right)\, \Aa_s\left(
\overline{x}\right)\right)\right)ds -\nabla (F\circ\chi_{s,t})\left(
\overline{x}\right)^{\prime}\, \Lambda_s\left(
\overline{x}\right)d\,\Ua_s
  \end{array}
 $$
 A proof of the above formula can be found in the articles~\cite{daprato-2,daprato-3}, see also~\cite{dp-ss}. Choosing the function $ F\left(
\overline{x}\right)=f(x)~z$, for some  twice differentiable function $f$ on $\RR^{m}$ with bounded derivatives and letting $z=1$ we check that
  $$
 \begin{array}{l}
 \displaystyle  d(f(X_{s,t}(x))Z_{s,t}(x))\\
 \\
 \displaystyle =-\left(\nabla (f(X_{s,t}(x))Z_{s,t}(x))^{\prime}~ a_s(x)+\frac{1}{2}~\tr\left(\nabla^2 (f(X_{s,t}(x))Z_{s,t}(x))~ \alpha_s(x)\right)\right)~ds\\
 \\
 \displaystyle\hskip3cm -~(f(X_{s,t}(x))Z_{s,t}(x))~b_s(x)^{\prime}\beta_s^{-1}~dY_s-\nabla (f(X_{s,t}(x))Z_{s,t}(x))^{\prime}~\sigma_s(x)~dW_s
  \end{array}
 $$
 This ends the proof of (\ref{backward-zakai-ps}). By localization arguments, the above result is also true for unbounded sensor functions.
 Integrating the flow of the signal we obtain (\ref{backward-zakai}). 
 This ends the proof of the theorem.
  \qed

We can also check (\ref{backward-zakai}) considering  a discrete time interval $[s,t]_h:=\{t_0,\ldots,t_{n-1}\}$ associated with some refining time mesh $t_{i+1}=t_i+h$ from $t_0=s$ to $t_{n}=t$, for some time step $h>0$. 
By (\ref{zakai}), for any $u\in [s,t]_h$  we compute $Q_{u,t}(f)$
from $Q_{u+h,t}(f)$ using the backward equation
\begin{eqnarray}
Q_{u,t}(f)&=&Q_{u+h,t}(f)+\left(Q_{u,u+h}-I\right)(Q_{u+h,t}(f))\nonumber\\
&\simeq&Q_{u+h,t}(f)+L_{u}(Q_{u+h,t}(f))\,h+Q_{u+h,t}(f)~b_{u}^{\prime}~\beta_{u}^{-1}\,(Y_{u+h}-Y_u)\label{ref-time-mesh}
\end{eqnarray}
  
For null sensor functions the evolution equation (\ref{backward-zakai-ps}) coincides with the backward It\^o formula discussed in~\cite{daprato-2,daprato-3,dp-ss,Kunita-82,kunita-b-82}. 

Choosing $f=1$ in (\ref{backward-zakai}) we recover the backward evolution of the likelihood function presented in~\cite{andersonRhodes,pardoux-3} (see formula (5.9) in~\cite{andersonRhodes} and equation (3.15) in~\cite{pardoux-3}). Arguing as in (\ref{perturbation-q-p}), using (\ref{backward-zakai}) we check the perturbation formulae given for any $s<t$ by,
$$
Q_{s,t}(f)=P_{s,t}(f)+\int_s^t~P_{s,u}\left(Q_{u,t}(f)~b_u^{\prime}\right)~\beta_u^{-1}\,dY_u
$$
Thus, for any $s<t$  the integral operator $Q_{s,t}(x_0,dx_1)$ has a density $(x_0,x_1)\mapsto q_{s,t}(x_0,x_1)$ given by  (\ref{perturbation-q-p}) and the integral formula,
\begin{equation}\label{perturbation-q-p-b}
q_{s,t}(x_0,x_1)=p_{s,t}(x_0,x_1)+\int_s^t~\left[\int~p_{s,u}(x_0,z)~q_{u,t}(z,x_1)~b_u^{\prime}(z)~dz\right]~\beta_u^{-1}\,dY_u
\end{equation}

  Using the same arguments as in the proof of Theorem~\ref{zakai-Q} we also have the following backward evolution equation.

\begin{prop}\label{zakai-Q-ov}
For any  twice differentiable function $f$ with bounded derivatives and for any
 $s\leq t$ we also have the backward equation
  \begin{eqnarray*}
  d\overline{\QQ}_{s,t}(f)(x) &=&-\left(\nabla\,\overline{\QQ}_{s,t}(f)(x)^{\prime}~ a_s(x)+\frac{1}{2}~\tr\left(
  \nabla^2\,\overline{\QQ}_{s,t}(f)(x)~ \alpha_s(x)\right)\right)~ds \\
&&\hskip1.7cm -\,\overline{\QQ}_{s,t}(f)(x)~\overline{b}_s(x)^{\prime}\beta_s^{-1}~d\overline{Y}_s-\nabla \overline{\QQ}_{s,t}(f)(x)^{\prime}~\sigma_s(x)~dW_s  
 \end{eqnarray*}
with terminal condition $\overline{\QQ}_{t,t}(f)=f$. In particular, we have the backward equation,
\begin{equation}\label{backward-zakai-2}
 d \overline{Q}_{s,t}(f)=-\left(L_s(\overline{Q}_{s,t}(f))\,ds+\overline{Q}_{s,t}(f)~ \overline{b}_s^{\prime}~\beta_s^{-1}\,d\overline{Y}_s\right)
\end{equation}
with terminal condition $\overline{Q}_{t,t}(f)=f$.
\end{prop}

Using the same notation as in (\ref{ref-time-mesh}), we also have the approximating backward equation
\begin{equation}\label{ref-time-over-Q}
\overline{Q}_{u,t}(f)\simeq\overline{Q}_{u+h,t}(f)+L_{u}(\overline{Q}_{u+h,t}(f))\,h+\overline{Q}_{u+h,t}(f)~\overline{b}_{u}^{\prime}~\beta_{u}^{-1}\,\left(\overline{Y}_{u+h}-\overline{Y}_u\right)
\end{equation}

\section{Smoothing semigroups and proof of the main result}\label{backward-sec}

This section is concerned with forward-backward evolution equations for the conditional smoothing distribution and the proof of our main result.

Let  $\Ka_{t,s}$  be the backward integral operator defined by,
 \begin{equation}\label{backward-Markov}
\Ka_{t,s}(f)(x):=\int~\pi_s(dz)~\frac{d\overline{Q}_{s,t}(z,\cdot)}{d\pi_t}(x)\,f(z)
 \end{equation}
For any $s\leq u\leq t$ we have the backward semigroup property,
 \begin{equation}\label{backward-semigroup-properties}
 \Ka_{t,s}=\Ka_{t,u}\,\Ka_{u,s}
 \end{equation}
which follows via,
\begin{eqnarray*}
(\Ka_{t,u}\,\Ka_{u,s})(f)(x)&=&\int~\pi_s(dx_0)~\overline{Q}_{s,u}(x_0,dx_1)~\frac{d\overline{Q}_{u,t}(x_1,\cdot)}{d\pi_t}(x)\,f(x_0)\\
&=&
\int~\pi_s(dx_0)~\frac{d\overline{Q}_{s,t}(x_0,\cdot)}{d\pi_t}(x)\,f(x_0)~=~\Ka_{t,s}(f)(x)
\end{eqnarray*}
and where we exploit the semigroup properties of the operators $\overline{Q}_{s,t}$. 

Also observe that for any $t>s>0$ the integral operator $\Ka_{t,s}(x_1,dx_0)$ has a density $(x_1,x_0)\mapsto k_{s,t}(x_1,x_0)$ w.r.t. the Lebesgue measure on $\RR^m$ given by,
$$
k_{t,s}(x_1,x_0):=p_s(x_0)~{\overline{q}_{s,t}(x_0,x_1)}/{p_t(x_1)}\quad \mbox{\rm with}\quad
\overline{q}_{s,t}(x_0,x_1)=q_{s,t}(x_0,x_1)/\pi_s(Q_{s,t}(1))
$$
The function $q_{s,t}$ denotes the density of the integral operator $Q_{s,t}$ discussed in (\ref{perturbation-q-p}) and (\ref{perturbation-q-p-b}).

Now, for any pair of functions $(f,g)$ we readily check the duality formula,
 \begin{equation}\label{duality}
\pi_s\left(f\,\overline{Q}_{s,t}(g)\right)=\pi_t\left(\Ka_{t,s}(f)\,g\right)
 \end{equation}

The following technical result is key in the proof of Theorem~\ref{theo-ref}.

 \begin{lem}\label{lem1}
  For any time parameter $s\leq t$ we have the forward-backward differential equation
\begin{equation}\label{fb-eq-2}
\partial_s\left(\pi_s\left(f~\overline{Q}_{s,t}(g)\right)\right)
=-
\pi_s\left(\overline{Q}_{s,t}(g)~\La_{s,\pi_s}(f)\right)
\end{equation}
with the second order differential operator
\begin{equation*}
 \La_{s,\pi_s}(f) \,:=\,-L_s(f)+\frac{1}{p_s}~~\sum_{1\leq i,j\leq m}\partial_{x_i}\left(p_s~\alpha^{i,j}_s~ \partial_{x_j}f\right)
\end{equation*}
\end{lem}

\proof
Observe that (\ref{fb-eq-2}) does not involve the derivatives of the function $g$.
Thus, up to a smooth mollifier's type approximation of the function $g$, it suffices to check (\ref{fb-eq-2}) for any pair of bounded and twice differentiable functions $f,g$ with bounded differentials. Arguing as in the proof of Proposition~\ref{prop-zakai} and Theorem~\ref{zakai-Q}, it suffices to prove the result for uniformly bounded sensor functions $b_u(x)$ on  $[s,t]\times\RR^m$, for any $s\leq t$.

In this situation, for any time horizon $t$, combining the  Kushner-Stratonovich equation (\ref{KS}) with 
the backward equation (\ref{backward-zakai-2}) for any $s\leq t$, we check the forward-backward evolution equation
\begin{equation}\label{fb-eq}
\partial_s\left(\pi_s\left(f~\overline{Q}_{s,t}(g)\right)\right)=\pi_s\left(L_s(f~\overline{Q}_{s,t}(g))-f~L_s(\overline{Q}_{s,t}(g))\right)
\end{equation}
The above equation can be proved using the backward It\^o-Ventzell formula in~\cite{dp-ss}. We use the same notation as in the proof of Theorem~\ref{zakai-Q}.  Let $\overline{Z}_{s,t}(x)$ be  the multiplicative functional  defined as $Z_{s,t}(x)$ by replacing  the function 
 $b_u$ and the observation It\^o-increment $dY_u$ by the centered function $ \overline{b}_u$ and the innovation 
  increment $d\overline{Y}_u$. 
  
Consider  the backward random field  $F_{s,t}$ with terminal condition $F_{t,t}\left(
\overline{x}\right)=f(x)g(x)z$ defined by the formula
$$
F_{s,t}\left(
\overline{x}\right):=f(x)~\overline{\QQ}_{s,t}(g)(x)~z\quad \mbox{\rm and we set}\quad
 \overline{\chi}_s:=\left(
 \begin{array}{cc}
 X_s\\
\overline{Z}_s
 \end{array}\right)\in (\RR^{m}\times\RR)\,.
$$
In this notation, we have
$$
\EE_0\left(F_{s,t}( \overline{\chi}_s)\right)=\EE_0\left(f(X_s)~\overline{Z}_s~\EE_0\left(\overline{\QQ}_{s,t}(g)(X_s)\vert~(X_s,Z_s)\right)\right)=\pi_s(f ~\overline{Q}_{s,t}(g))
$$
Observe that
$$
F_{s,t}\left(
\overline{x}\right)
 =f(x)~(F\circ \overline{\chi}_{s,t})\left(
\overline{x}\right)
 $$
 with the function
 $$
F\left(
\overline{x}\right):=g(x)~z\quad\mbox{\rm and the stochastic flow}\quad
\overline{\chi}_{s,t}(x,z):=\left(\begin{array}{c}
 X_{s,t}(x)\\
\overline{Z}_{s,t}(x)~z
 \end{array}\right)
$$
Following the proof of Theorem~\ref{zakai-Q}, we check that
$$
dF_{s,t}\left(
\overline{x}\right)=f(x)~d(F\circ \overline{\chi}_{s,t})\left(
\overline{x}\right)=-\left(\Ga_{s,t}\left(
\overline{x}\right)~ds+\Ha_{s,t}\left(
\overline{x}\right)~d\,\Ua_s\right)
$$
with the drift function
$$
\Ga_{s,t}\left(
\overline{x}\right):=~f(x)~z~\left(\nabla\, \overline{\QQ}_{s,t}(g)(x)^{\prime}~ a_s(x)+\frac{1}{2}~\tr\left(
 \nabla^2 \overline{\QQ}_{s,t}(g)(x)^{\prime}~ \alpha_s(x)\right)\right)
$$
and the diffusion term
$$
\Ha_{s,t}\left(
\overline{x}\right)~d\,\Ua_s:=~f(x)~z~\left(
 \nabla\,\overline{\QQ}_{s,t}(g)(x)^{\prime}~\sigma_s(x)~dW_s+\overline{\QQ}_{s,t}(g)(x)~b_s(x)^{\prime}\beta_s^{-1}~dY_s
 \right)
$$
Applying the backward It\^o-Ventzell formula~\cite{dp-ss} we check that
$$
 dF_{s,t}( \overline{\chi}_s)=(dF_{s,t})( \overline{\chi}_s)+\nabla F_{s,t}( \overline{\chi}_s)^{\prime}~d\chi_s+\frac{1}{2}~\tr\left(\nabla^2 F_{s,t}(\chi_s)^{\prime}~\Aa_t
\left( \overline{\chi}_s\right)\right)ds
$$
from which we conclude that
 $$
 \begin{array}{l}
\displaystyle dF_{s,t}( \overline{\chi}_s)=\overline{Z}_s\left(\nabla\, (\overline{\QQ}_{s,t}(g)(x)~f(x))^{\prime}_{\vert x=X_s}-f(X_s)~\overline{Z}_s~ \nabla\, \overline{\QQ}_{s,t}(g)(X_s)^{\prime}\right)\sigma_s(X_s)\,dW_s\\
\\
\displaystyle~~-f(X_s)\overline{Z}_s~\left( \nabla\,\overline{\QQ}_{s,t}(g)(X_s)^{\prime}~ a_s(X_s)~ds+\frac{1}{2}~ \tr\left(
\nabla^2\, \overline{\QQ}_{s,t}(g)(X_s)~ \alpha_s(X_s)\right)\right)ds\\
\\
\displaystyle~~+\overline{Z}_s\left(\nabla\, (\overline{\QQ}_{s,t}(g)(x)~f(x))^{\prime}_{\vert x=X_s}~a_s(X_s)~ds+~\frac{1}{2}~ \tr\left(\nabla^2\, (\overline{\QQ}_{s,t}(g)(x)~f(x))^{\prime}_{\vert x=X_s}~ \alpha_s(X_s)\right)\right)ds
  \end{array}
$$
We end the proof of (\ref{fb-eq}) by simple integration.

To take the final step, we recall the integration by parts formula
 \begin{equation*}
  L_t(fg)=f~L_t(g)+g ~L_t(f)+\Gamma_{L_t}(f,g)
 \end{equation*}
with the carr\'e-du-champ (a.k.a. square field) operator $   \Gamma_{L_t}$ associated with the generator $L_t$ defined by
 \begin{equation*}
   \Gamma_{L_t}(f,g):=(\nabla f)^{\prime}\alpha_t\nabla g
 \end{equation*}
  Combining (\ref{fb-eq}) with the above formula we check that
\begin{equation*}
\partial_s\left(\pi_s\left(f~\overline{Q}_{s,t}(g)\right)\right)=\pi_s\left(L_s(f)~\overline{Q}_{s,t}(g)\right)+\pi_s\left(
   \Gamma_{L_s}\left(\overline{Q}_{s,t}(g),f\right)\right)
\end{equation*}
On the other hand, by an integration by parts we have
\begin{equation*}
\pi_s\left(
   \Gamma_{L_s}\left(\overline{Q}_{s,t}(g),f\right)\right)
=-\sum_{i,j}\int~p_s(x)~
 \overline{Q}_{s,t}(g)(x) \frac{1}{p_s(x)} \partial_{x_i}~\left(p_s(x)~\alpha_t^{i,j}~\partial_{x_j} f(x)\right)dx
\end{equation*}
This ends the proof of the lemma.
\qed

Another approach for finding (\ref{fb-eq}) is to use for any $u\in [s,t]_h$ the decomposition
 \begin{equation}
\begin{array}{l}
\pi_{u+h}\left(f\,\overline{Q}_{u+h,t}(g)\right)-\pi_u\left(f\,\overline{Q}_{u,t}(g)\right)\\
\\
\qquad =\,\pi_{u}\left(f\left(\overline{Q}_{u+h,t}-\overline{Q}_{u,t}\right)(g)\right)+(\pi_{u+h}-\pi_u)\left(f~\overline{Q}_{u+h,t}(g)\right)
\end{array}
 \end{equation}
 
Note that $\pi_{u}$ depends on the observations $(Y_s-Y_0)$ from $s=0$ up to time $s=u$, while the increment
$\overline{Q}_{u,t}$ is computed backward in time and only depends on the observations $(Y_s-Y_u)$ from $s>u$ up to $s=t$. Conversely, 
$\pi_{u+h}$ depends on the observations $(Y_s-Y_0)$ from $s=0$ up to time $s=u+h$, while $\overline{Q}_{u+h,t}$ 
is computed backward in time and only depends on the observations $(Y_s-Y_{u+h})$ from $s>u+h$, up to time $s=t$.

Following the  two-sided stochastic integration calculus developped by Pardoux and Protter in~\cite{pardoux-protter} (see also~\cite{dp-ss} for extended versions to interpolating stochastic flows),
combining the forward (\ref{ref-time-over-pi}) with the backward equation (\ref{ref-time-over-Q}), when $h\simeq 0$  we can check the approximation,
$$
\displaystyle\sum_{u\in [s,t]_h}\left\{\pi_{u+h}\left(f~\overline{Q}_{u+h,t}(g)\right)-\pi_u\left(f~\overline{Q}_{u,t}(g)\right)-
 \pi_u\left(L_u(f~\overline{Q}_{u+h,t}(g))-f~L_{u}(\overline{Q}_{u+h,t}(g))\right)\, h\right\}\simeq ~0
$$

\subsection{Proof of Theorem~\ref{theo-ref}}\label{theo-ref-proof}

With the definition in (\ref{backward-Markov}) we have,
 \begin{equation}\label{backward-semigroup-smoother}
 \pi_{t,s}(dx)=(\pi_t\,\Ka_{t,s})(dx)=\pi_s(dx)\,\overline{Q}_{s,t}(1)(x)
 \end{equation}
The  formulation of the conditional distribution $\pi_{t,s}$ of $X_s$ given $\Ya_t$ in (\ref{backward-semigroup-smoother}) is rather well known, see e.g. Theorem 3.7 and Corollary 3.8 in~\cite{pardoux}, as well as equation (3.9) in~\cite{andersonRhodes}. The proof of this formula is a direct consequence of (\ref{backward-Markov}). With (\ref{backward-semigroup-properties}) we have,
$$
  \pi_t \Ka_{t,s}=\pi_{t,u}\Ka_{u,s}=\pi_{t,s} 
$$
Thus with $\Ka_{t,s}$ as defined in (\ref{backward-Markov}) we immediately have the transport equation in (\ref{theo-eq-pi-K}).

It remains to show that this integral operator (as defined in (\ref{backward-Markov})) is also the Markov transition kernel of the backward diffusion flow in (\ref{theo-eq}). The rest of the proof of Theorem~\ref{theo-ref} is a consequence of the duality formula (\ref{duality}) and Lemma~\ref{lem1}.

Rewritten in a slightly different form, the duality formula (\ref{duality}) reads as follows,
$$
\EE\left(f(X_s)~g(X_t)~|~\Ya_t\right)=\EE\left(\Ka_{t,s}(f)(X_t)~g(X_t)~|~\Ya_t\right)
$$
 This implies that
$$
 \Ka_{t,s}(f)(X_t)=\EE\left(f(X_s)~|~X_t,~\Ya_t\right)
$$
Finally, combining (\ref{fb-eq-2}) with the duality formula (\ref{duality}) we have
\begin{equation*}
\pi_t\left(g~\partial_s\Ka_{t,s}(f)\right)=-\pi_t\left(g~\Ka_{t,s}(\La_{s,\pi_s}(f))\right)
\end{equation*}
Since the above formula is valid for any test function $g$ and $\pi_t$ has a bounded positive density, we check the backward Kolmogorov equation
\begin{equation}\label{ref-end-2}
\partial_s\Ka_{t,s}(f)(x)=
-\Ka_{t,s}(\La_{s,\pi_s}(f))(x)
\end{equation}
with terminal condition $\Ka_{t,t}(f)=f$, when $s=t$, for almost every $x\in\RR^m$ (and almost surely w.r.t. the law of the observation process from the origin up to the time $t$). Since both terms in (\ref{ref-end-2}) are continuous, the above equality holds for any $x\in\RR^m$, almost surely. 

We now complete the proof by showing that the integral operator $\Ka_{t,s}(x,dz)$ (defined in (\ref{backward-Markov})) does indeed coincide with the transition kernel associated with the flow $\Xa_{t,s}(x)$ in (\ref{theo-eq}). Firstly, observe that (\ref{ref-end-2}) coincides with the backward Kolmogorov equation (\ref{fb-eq-leondes}) associated with the transition semigroup of the stochastic flow $\Xa_{t,s}(x)$. Denote this transition semigroup by $\overline{\Ka}_{t,s}(x,dz)$ temporarily.

By the semigroup properties of $\overline{\Ka}_{t,s}$, for any $s\leq u\leq t$ and any smooth function $f$ we have
$$
\partial_u\overline{\Ka}_{t,s}(f)=0=\partial_u(\overline{\Ka}_{t,u}(\overline{\Ka}_{u,s}(f)))=-\overline{\Ka}_{t,u}(\La_{u,\pi_u}(\overline{\Ka}_{u,s}f))+\overline{\Ka}_{t,u}(\partial_u\overline{\Ka}_{u,s}(f))
$$
Choosing $u=t$ we obtain the forward equation
$$
\partial_t\overline{\Ka}_{t,s}(f)=
\La_{t,\pi_t}(\overline{\Ka}_{t,s}(f))
$$
Arguing as above, this implies that
$$
\partial_u(\Ka_{t,u}(\overline{\Ka}_{u,s}(f)))=-\Ka_{t,u}(\La_{u,\pi_u}(\overline{\Ka}_{u,s}f))+\Ka_{t,u}(\La_{u,\pi_u}(\overline{\Ka}_{u,s}(f)))=0
$$
Integrating over the interval $[s,t]$ we check that $\Ka_{t,s}=\overline{\Ka}_{t,s}$.
This ends the proof of Theorem~\ref{theo-ref}.\qed

\end{document}